\documentclass[preprint,12pt]{elsarticle}

\usepackage{amssymb}
\usepackage{hyperref}
\usepackage{amsmath}
\usepackage{fixltx2e}
\usepackage{footnote}
\usepackage{multicol}
\usepackage{booktabs}
\usepackage{graphicx}
\pdfoutput=1
\usepackage[flushleft]{threeparttable}
\usepackage[font=small,labelfont=bf]{caption}
\newtheorem{thm}{Theorem}
\newtheorem{lemma}[thm]{Lemma}
\newdefinition{rmk}{Remark}
\newproof{proof}{Proof}

\journal{XXX}

\begin{document}

\begin{frontmatter}


 \title{Further Results on Topological Indices of Two Corona Variants of Graph Products}
%
%
%

%
\author[label1]{Rao Tabish}
\ead{ 03338431034tabish@gmail.com }
\author[label1]{Muhammad Kamran Jamil}
\ead{m.kamran.sms@gmail.com}

\author[label2]{Muhammad Waheed\fnref{label5}}
\ead{waheed.i35401@gmail.com}

%
\fntext[label5]{Corresponding author}
\address[label1]{Department of Mathematics, Riphah Institute of Computing and Applied Sciences, Riphah International University, Lahore, Pakistan.}
\address[label2]{Govt Higher Secondary School Nangal Sahdan Muridke, District Sheikhupura, Pakistan.}

\begin{abstract}
	
 Graph operations are crucial for building complicated network architectures from simple graphs. In \cite{bjwa}, authors defined two new variants of Corona product and investigated their some topological indices. In this paper, we extended the work and found the formulas of forgotten, first hyper Zagreb, and reduced second Zagreb indices for Corona join product and sub-division vertex join products of graphs.

\end{abstract}

\begin{keyword}
graph operations, topological indices, subdivision, join, corona product



\end{keyword}

\end{frontmatter}


\section{Introduction}

In Mathematics \cite{1}-\cite{2}, graph theory is the study of graphs which are mathematical structures used to model pairwise connection between objects. Graph theory is applied in the various fields like
software engineering, biology, chemistry and operation
research .

Let $G$ be a simple, connected graph with vertex set $V(G)$ and edge set $E(G)$. In graph $G$, the number of vertices and number of edges are called order $n$ and size $m$, respectively. A graph of order $n$ and size $m$ will be denoted by $G(n,m)$. Suppose for any vertex $v\in V(G)$, the degree of vertex $v$ is the number of edges incident on the vertex $v$ and it is written as $d_{G}(v)$ or simply $d(v)$. For a simple graph $G$, the sub-division of the graph $G$ is denoted by $S(G)$ and obtained by inserting a new vertex on every edge of $G$. 

Topological index is a numeric value which is associated with a chemical structure of a certain chemical compound. Topological index also called a molecular structure  descriptor or graph
theoretical descriptor \cite{3}. The first and second Zagreb indices are defined as:
\begin{align*}
	M_1(G)=\sum_{v \in V(G)}d_{G}(v)^2=\sum_{uv\in E(G)}[d_{G}(u)+d_{G}(v)]
\end{align*}
\begin{align*}
	M_2(G)=\sum_{uv \in E(G)}{d_G(u)d_G(v)}
\end{align*}
In 1972 \cite{4}, these topological indices applied for the first time to find the total $\pi$ energy of molecular graphs. Later, the Zagreb indices developed \\important applications in QSPR/QSAR studies and a lot of research paper have been published on these, see \cite{5}-\cite{10}. 

The forgotten topological index or F-index \cite{11,12} is defined as 

\[F(G)=\sum_{v\in V(G)}d_{G}(v)^3=\sum_{uv\in E(G)}[d_{G}(u)^2+d_{G}(v)^2]\]

In 2013 Shirdel and co-authors \cite{13} investigated a new version of Zagreb indices called first hyper Zagreb index and it is defined as

\[HM_{1}(G)=\sum_{uv\in E(G)}[d_{G}(u)+d_{G}(v)]^2\]

In 2019, Abdussakir \cite{14} proposed a new extension of second Zagreb index which is defined as

\[RM_2(G)=\sum_{uv \in E(G)}[d_G(u)-1][d_G(v)-1]\]

A large number of the topological indices have discovered applications as intends to demonstrate substance, drug and different properties of atoms.\\

The graph operation, especially graph products, plays significant role not only in pure
and applied mathematics, but also in computer science. For instance, the cartesian product provide a significant  model for connecting computers \cite{15}.

Let $G_1(n_1,m_1)$ and $G_2(n_2,m_2)$, be two connected simple graphs. Corona product is shown by $G_1\circ G_2$ and it is obtained by taking one copy of $G_{1}$ and $n_{1}$ copies of $G_{2}$, and joining the each vertex of $i$-$th$ copy of $G_{2}$ with the vertex $v_{i}\in n_1$.

Suppose $G_{1}$ and $G_{2}$ be two graphs on disjoint vertex sets. The graph $G_{1}+G_{2}$ has vertex set $V_{1}\cup V_{2}$ and the edge set $E(G_{1}+G_{2})=E_{1}\cup E_{2}\cup \{uv: u\in V_{1}, v\in V_{2}\}$. The sum of two graphs is also known as join \cite{16}.\\
P. Sarkar in \cite{17} computed the forgotten topological index of double corona of graphs related to the different subdivision graphs. Authors in \cite{18}, investigated the upper bounds on the multiplicative Zagreb indices of some product of graphs. Nilanjan De \cite{19}  presented the f co-index of some graph operations. Jamil et al. \cite{20}, found the exact formulas of the first reformulated Zagreb index for Cartesian product, composition, join, corona product, splice, link and chain of graphs. Some graph operations and their topological indices are presented in \cite{21}-\cite{27}.\\
  Now we defined variants of these graphs product.

\begin{figure}[h]
	\centering
	\includegraphics[width=0.4\linewidth]{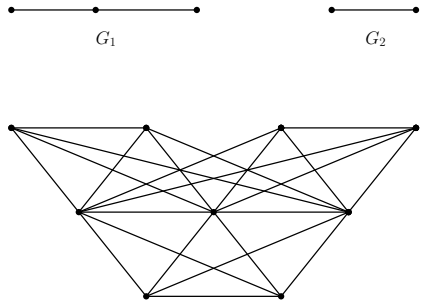}
	\caption{Corona join product $G_{1}$$\oplus G_{2}$}
	\label{fig:corona-join-product}
\end{figure}

\textbf{Corona join product:}
Let we have two simple connected graphs\\ $G_{1}=(n_{1},m_{1})$ and $G_{2}=(n_{2},m_{2})$, then corona join graph of $G_1$ and $G_2$ is denoted by $G_{1}$$\oplus G_{2}$ and obtained by taking one copy of $G_1$, $n_{1}$ copies of $G_2$ and by joining each vertex of the $i-th$ copy of $G_{2}$ with all vertices of $G_{1}$. Corona join product of two graphs is presented in the Figure \ref{fig:corona-join-product}.

\begin{figure}[h]
	\centering
	\includegraphics[width=0.4\linewidth]{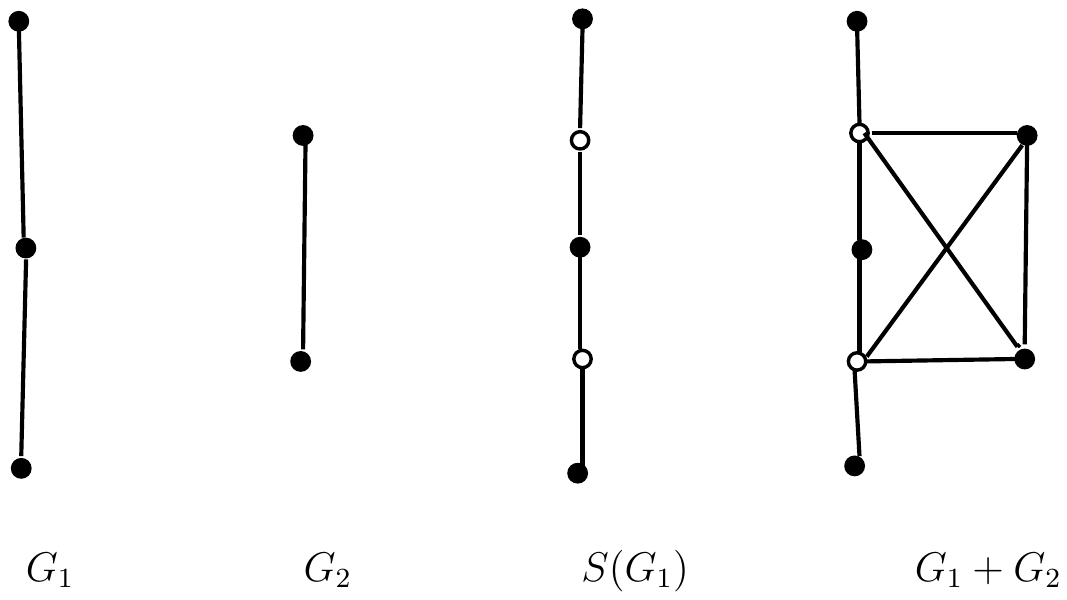}
	\caption{Subdivision-vertex join graph $G_1\dotplus G_2$}
	\label{fig:subdivision-vertex-join-graph}
\end{figure}

\textbf{Sub-division vertex join graph:}
Let we have three simple connected graphs $G_{1}=(n_{1},m_{1})$, $G_{2}=(n_{2},m_{2})$ and $S(G_{1})=(n_{1}^{\prime},m_{1}^{\prime})$, then sub-division vertex join graph is denoted by $G_1\dotplus G_2$ and obtained by joining the each new vertex of $S(G_1)$ to all vertices of $G_2$. Sub-division vertex join graph is presented in the Figure \ref{fig:subdivision-vertex-join-graph}.

\section{Main Results}
In this section, we present the main results. The following lemma's are useful to obtain the exact expressions of topological  indices of new variants of graph products. The proofs of the following two lemmas are directly from the definitions of Corona join product $G_{1}$$\oplus G_{2}$ and sub-division vertex join $G_{1}$$\dotplus G_{2}$.     
\begin{lemma} \label{cjp}
Let we have two simple connected graphs $G_{1}=(n_{1},m_{1})$ and $G_{2}=(n_{2},m_{2})$, then the degree behavior of vertices in the graph $G_{1}$$\oplus G_{2}$ is 
\[
d_{G_{1}\oplus G_{2}}(v)=
\begin{cases}
	d_{G_{1}}(v)+n_{1}n_{2},  &\text{if $v\in V(G_1)$;}\\
	d_{G_{2}}(v)+n_{1},  &\text{if $v\in V(G_2)$;}\\  
\end{cases}
\]  	
\end{lemma}

\begin{lemma} \label{svjg}
	Let we have three simple connected graphs $G_{1}=(n_{1},m_{1})$, \\$G_{2}=(n_{2},m_{2})$ and $S(G_{1})=(n_{1}^{\prime},m_{1}^{\prime})$, then the degree behavior of vertices in the graph $G_{1}\dotplus G_{2}$ is 
	\[
	d_{G_{1}\dotplus G_{2}}(v)=
	\begin{cases}
		d_{G_{1}}(v),  &\text{if $v\in V(G_1)$;}\\
		2+n_2,  &\text{if $v\in V_{s}(G_1)$;}\\ 
		d_{G_{2}}(v)+m_{1},  &\text{if $v\in V(G_2)$;}\\  
	\end{cases}
	\]  	
\end{lemma}

\begin{thm} \label{fi1}
Let we have two simple connected graphs $G_{1}=(n_{1},m_{1})$ and $G_{2}=(n_{2},m_{2})$, then the forgotten topological index of Corona join product $G_{1}\oplus G_{2}$ is given as
\begin{align*} 
	F(G_{1}\oplus G_{2})=&F(G_1)+n_{1}F(G_2)+3n_{1}n_{2}M_{1}(G_{1})+3n_{1}^2M_{1}(G_2)+6m_{1}n_{1}^2n_{2}^2\\&+6m_{2}n_{1}^3+n_{1}^4n_{2}(n_{2}^2+1)
\end{align*} 

\begin{proof}
	
   From the definition of the forgotten topological index we have
   \[F\Big(G_{1}\oplus G_{2}\Big)=\sum_{v\in V(G_{1}\oplus G_{2})}d_{G_{1}\oplus G_{2}}(v)^3\]
   Now we apply the Lemma \ref{cjp}
   \begin{align*}
   	&=\sum_{v\in V(G_{1})}\Big(d_{G_{1}}(v)+n_{1} 
   	n_{2}\Big)^3+\sum_{v\in V(G_{1})}\sum_{v\in V(G_{2})}\Big(d_{G_{2}}(v)+n_{1}\Big)^3\\
   	&=\sum_{v\in V(G_{1})}\Big(d_{G_{1}}(v)^3+(n_{1}n_{2})^3+3d_{G_{1}}(v)^2(n_{1} 
   	n_{2})+3d_{G_{1}}(v)(n_{1}n_{2})^2\Big)+\\&\sum_{v\in V(G_{1})}\sum_{v\in V(G_{2})}\Big(d_{G_{2}}(v)^3+n_{1}^3+3d_{G_{2}}(v)^2(n_{1})+3d_{G_{2}}(v)(n_{1})^2\Big)\\
   	&=\sum_{v\in V(G_{1})}d_{G_{1}}(v)^3+(n_{1}n_{2})^3\sum_{v\in V(G_{1})}1+3(n_{1} 
   	n_{2})\sum_{v\in V(G_{1})}d_{G_{1}}(v)^2+3(n_{1}n_{2})^2\\&\sum_{v\in V(G_{1})}d_{G_{1}}(v)+\sum_{v\in V(G_{1})}1\sum_{v\in V(G_{2})}d_{G_{2}}(v)^3+n_{1}^3\sum_{v\in V(G_{1})}1\sum_{v\in V(G_{2})}1+3n_{1}\sum_{v\in V(G_{1})}1\\&\sum_{v\in V(G_{2})}d_{G_{2}}(v)^2+3n_{1}^2\sum_{v\in V(G_{1})}1\sum_{v\in V(G_{2})}d_{G_{2}}(v)\\
   	&=F(G_1)+n_{1}^4n_{2}^3+3n_{1}n_{2}M_{1}(G_{1})+6m_{1}n_{1}^2n_{2}^2+n_{1}F(G_2)+n_{1}^4n_{2}+3n_{1}^2M_{1}(G_2)\\&+6m_{2}n_{1}^3\\
   	&=F(G_1)+n_{1}F(G_2)+3n_{1}n_{2}M_{1}(G_{1})+3n_{1}^2M_{1}(G_2)+6m_{1}n_{1}^2n_{2}^2+6m_{2}n_{1}^3\\&+n_{1}^4n_{2}(n_{2}^2+1)
   \end{align*}
which is our required result.
   \end{proof}	
\end{thm}

Here is the example. By using the statement of theorem \ref{fi1}, we get
\begin{align*}
	F(P_{l}\oplus C_{m})=&8l-14+8lm+3lm(4l-6)+12ml^2+6l^2m^2(l-1)\\&+6ml^3+l^4m(m^2+1)
\end{align*}

\begin{thm} \label{fi2}
	Let we have three simple connected graphs $G_{1}=(n_{1},m_{1})$, \\$G_{2}=(n_{2},m_{2})$ and $S(G_{1})=(n_{1}^{\prime},m_{1}^{\prime})$, then the forgotten topological index of sub-division vertex join $G_{1}\dotplus G_{2}$ is given as
	\begin{align*} 
		F(G_{1}\dotplus G_{2})=&F(G_1)+F(G_2)+3m_{1}M_{1}(G_{2})+m_{1}(2+n_{2})^3+6m_{2}m_{1}^2+m_{1}^3n_{2}
	\end{align*}
	
	\begin{proof}
		From the definition of the forgotten topological index we have
		\[F\Big(G_{1}\dotplus G_{2}\Big)=\sum_{v\in V(G_{1}\dotplus G_{2})}d_{G_{1}\dotplus G_{2}}(v)^3\]
		Now we apply the Lemma \ref{svjg}
		\begin{align*}
			&=\sum_{v\in V(G_{1})}d_{G_1}(v)^3+\sum_{v\in V_{s}(G_{1})}\Big(2+n_{2}\Big)^3 +\sum_{v\in V(G_{2})}\Big(d_{G_{2}}(v)+m_{1}\Big)^3\\
			&=\sum_{v\in V(G_{1})}d_{G_1}(v)^3+\sum_{v\in V_{s}(G_{1})}\Big(2+n_{2}\Big)^3+\sum_{v\in V(G_{2})}\Big(d_{G_{2}}(v)^3+m_{1}^3+3d_{G_{2}}(v)^2m_{1}\\&+3d_{G_{2}}(v)m_{1}^2\Big)\\
			&=\sum_{v\in V(G_{1})}d_{G_1}(v)^3+\Big(2+n_{2}\Big)^3\sum_{v\in V_{s}(G_{1})}1+\sum_{v\in V(G_{2})}d_{G_{2}}(v)^3+m_{1}^3\sum_{v\in V(G_{2})}1+3m_{1}\\&\sum_{v\in V(G_{2})}d_{G_{2}}(v)^2+3m_{1}^2\sum_{v\in V(G_{2})}d_{G_{2}}(v)\\
			&=F(G_1)+m_{1}(2+n_{2})^3+F(G_2)+m_{1}^3n_{2}+3m_{1}M_{1}(G_{2})+6m_{1}^2m_{2}\\
			&=F(G_1)+F(G_2)+3m_{1}M_{1}(G_{2})+m_{1}(2+n_{2})^3+6m_{2}m_{1}^2+m_{1}^3n_{2}
		\end{align*}
		which is our required result.
	\end{proof}
\end{thm}

Here is the example. By using the statement of theorem \ref{fi2}, we get
\begin{align*}
	F(P_{l}\dotplus C_{m})=&8l-14+8m+12m(l-1)+(l-1)(2+m)^3\\&+6m(l-1)^2+m(l-1)^3
\end{align*}

\begin{thm} \label{H1}
	Let we have two simple connected graphs $G_{1}=(n_{1},m_{1})$ and $G_{2}=(n_{2},m_{2})$, then the first hyper Zagreb index of Corona join product $G_{1}\oplus G_{2}$ is given as
	\begin{align*}
		&HM_{1}(G_{1}\oplus G_{2})=HM_{1}(G_1)+n_{1}HM_{1}(G_2)+5n_{1}\Big[n_{2}M_{1}(G_1)+n_{1}M_{1}(G_2)\Big]\\&+4n_{1}\Big[m_{1}n_{1}n_{2}^2+m_{2}n_{1}^2+2m_{1}m_{2}\Big]+n_{1}^2(n_{2}+1)\Big[4m_{1}n_{2}+4m_{2}n_{1}+n_{1}^2n_{2}(n_{2}+1)\Big]
	\end{align*}
	
	\begin{proof}
		From the definition of the first hyper Zagreb index we have
		\begin{align*}
			HM_{1}\Big(G_{1}\oplus G_{2}\Big)=&\sum_{uv\in E(G_{1}\oplus G_{2})}\Big[d_{G_{1}\oplus G_{2}}(u)+d_{G_{1}\oplus G_{2}}(v)\Big]^2
		\end{align*}
		Now we apply the Lemma \ref{cjp}
		\begin{align*}
			&=\sum_{uv\in E(G_{1})}\Big[d_{G_{1}}(u)+n_{1}n_{2}+d_{G_{1}}(v)+n_{1}n_{2}\Big]^2+n_{1}\sum_{uv\in E(G_{2})}\Big[d_{G_{2}}(u)+n_{1}\\&+d_{G_{2}}(v)+n_{1}\Big]^2+n_{1}\sum_{\substack{uv\in E(G_{1}\oplus G_{2})\\ u\in V(G_1)\\ v\in V(G_2)}}\Big[d_{G_{1}}(u)+n_{1}n_{2}+d_{G_{2}}(v)+n_{1}\Big]^2\\
			&=\sum_{uv\in E(G_{1})}\Big[\Big(d_{G_{1}}(u)+d_{G_{1}}(v)\Big)+2n_{1}n_{2}\Big]^2+n_{1}\sum_{uv\in E(G_{2})}\Big[\Big(d_{G_{2}}(u)+d_{G_{2}}(v)\Big)+2n_{1}\Big]^2\\&+n_{1}\sum_{\substack{uv\in E(G_{1}\oplus G_{2})\\ u\in V(G_1)\\ v\in V(G_2)}}\Big[\Big(d_{G_{1}}(u)+d_{G_{2}}(v)\Big)+\Big(n_{1}(n_{2}+1)\Big)\Big]^2\\
			&=\sum_{uv\in E(G_{1})}\Big[\Big(d_{G_{1}}(u)+d_{G_{1}}(v)\Big)^2+4n_{1}^2n_{2}^2+4n_{1}n_{2}\Big(d_{G_{1}}(u)+d_{G_{1}}(v)\Big)\Big]\\&+n_{1}\sum_{uv\in E(G_{2})}\Big[\Big(d_{G_{2}}(u)+d_{G_{2}}(v)\Big)^2+4n_{1}^2+4n_{1}\Big(d_{G_{2}}(u)+d_{G_{2}}(v)\Big)\Big]\\&+n_{1}\sum_{\substack{uv\in E(G_{1}\oplus G_{2})\\ u\in V(G_1)\\ v\in V(G_2)}}\Big[d_{G_{1}}(u)^2+d_{G_{2}}(v)^2+2d_{G_{1}}(u)d_{G_{2}}(v)+n_{1}^2(n_{2}+1)^2+2n_{1}(n_{2}+1)\\&d_{G_{1}}(u)+2n_{1}(n_{2}+1)d_{G_{2}}(v)\Big]\\
			&=\sum_{uv\in E(G_{1})}\Big(d_{G_{1}}(u)+d_{G_{1}}(v)\Big)^2+4n_{1}^2n_{2}^2\sum_{uv\in E(G_{1})}1+4n_{1}n_{2}\sum_{uv\in E(G_{1})}\Big(d_{G_{1}}(u)+d_{G_{1}}(v)\Big)\\&+n_{1}\sum_{uv\in E(G_{2})}\Big(d_{G_{2}}(u)+d_{G_{2}}(v)\Big)^2+4n_{1}^3\sum_{uv\in E(G_{2})}1+4n_{1}^2\sum_{uv\in E(G_{2})}\Big(d_{G_{2}}(u)+d_{G_{2}}(v)\Big)\\&+n_{1}\sum_{u\in V(G_1)}d_{G_1}(u)^2\sum_{v\in V(G_2)}1+n_{1}\sum_{v\in V(G_2)}d_{G_2}(v)^2\sum_{u\in V(G_1)}1+2n_{1}\sum_{u\in V(G_1)}d_{G_1}(u)\\&\sum_{v\in V(G_2)}d_{G_2}(v)+n_{1}^3(n_{2}+1)^2\sum_{u\in V(G_1)}1\sum_{v\in V(G_2)}1+2n_{1}^2(n_{2}+1)\sum_{u\in V(G_1)}d_{G_1}(u)\sum_{v\in V(G_2)}1\\&+2n_{1}^2(n_{2}+1)\sum_{v\in V(G_2)}d_{G_2}(v)\sum_{u\in V(G_1)}1\\
		\end{align*}
		\begin{align*}
		&=HM_{1}(G_{1})+4n_{1}^2n_{2}^2m_{1}+4n_{1}n_{2}M_{1}(G_{1})+n_{1}HM_{1}(G_{2})+4n_{1}^3m_{2}+4n_{1}^2M_{1}(G_{2})\\&+n_{1}n_{2}M_{1}(G_{1})+n_{1}^2M_{1}(G_{2})+8m_{1}m_{2}n_{1}+n_{1}^4n_{2}(n_{2}+1)^2+4m_{1}n_{1}^2n_{2}(n_{2}+1)\\&+4m_{2}n_{1}^3(n_{2}+1)\\
			&=HM_{1}(G_1)+n_{1}HM_{1}(G_2)+5n_{1}\Big[n_{2}M_{1}(G_1)+n_{1}M_{1}(G_2)\Big]+4n_{1}\\&\Big[m_{1}n_{1}n_{2}^2+m_{2}n_{1}^2+2m_{1}m_{2}\Big]+n_{1}^2(n_{2}+1)\Big[4m_{1}n_{2}+4m_{2}n_{1}+n_{1}^2n_{2}(n_{2}+1)\Big]
		\end{align*}
		
		which is our required result.
	\end{proof}	
\end{thm}

Here is the example. By using the statement of theorem \ref{H1}, we get
\begin{align*}
	&HM_{1}(P_{l}\oplus C_{m})=16l-30+16lm+5l\Big[m(4l-6)+4lm\Big]+\\&4l\Big[lm^2(l-1)+ml^2+2m(l-1)\Big]+l^2(m+1)\Big[4m(l-1)+4ml+l^2m(m+1)\Big]
\end{align*}

\begin{thm} \label{H2}
	Let we have three simple connected graphs $G_{1}=(n_{1},m_{1})$, \\$G_{2}=(n_{2},m_{2})$ and $S(G_{1})=(n_{1}^{\prime},m_{1}^{\prime})$, then the first hyper Zagreb index of sub-division vertex join $G_{1}\dotplus G_{2}$ is given as
	\begin{align*}
		&HM_{1}\Big(G_{1}\dotplus G_{2}\Big)=m_{1}\Big[5M_{1}(G_2)+4m_{1}m_{2}\Big]+HM_{1}(G_2)+m_{1}^{\prime}\Big(2+n_{2}\Big)^2+m_{1}\\&\Big(m_{1}+n_{2}+2\Big)\Big[n_{2}\Big(m_{1}+n_{2}+2\Big)+4m_{2}\Big]+\sum_{uv\in E(S(G_1))}\Big[d_{G_1}(u)^2+\Big(4+2n_{2}\Big)d_{G_1}(u)\Big]
	\end{align*}

	\begin{proof}
		From the definition of the first hyper Zagreb index we have
		\begin{align*}
			HM_{1}\Big(G_{1}\dotplus G_{2}\Big)=&\sum_{uv\in E(G_{1}\dotplus G_{2})}\Big[d_{G_{1}\dotplus G_{2}}(u)+d_{G_{1}\dotplus G_{2}}(v)\Big]^2
		\end{align*}
		Now we apply the Lemma \ref{svjg}
		\begin{align*}
			&=\sum_{\substack{uv\in E(S(G_1))\\ u\in V(G_1)\\ v\in V_{s}(G_1)}}\Big[d_{G_1}(u)+\Big(2+n_{2}\Big)\Big]^2+\sum_{uv\in E(G_2)}\Big[\Big(d_{G_2}(u)+d_{G_2}(v)\Big)+2m_{1}\Big]^2\\&+\sum_{\substack{uv\in E(G_1\dotplus G_2)\\ u\in V_{s}(G_1)\\ v\in V(G_2)}}\Big[d_{G_2}(v)+\Big(m_{1}+n_{2}+2\Big)\Big]^2\\
		\end{align*}
		\begin{align*}
			&=\sum_{\substack{uv\in E(S(G_1))\\ u\in V(G_1)\\ v\in V_{s}(G_1)}}\Big[d_{G_1}(u)^2+\Big(2+n_{2}\Big)^2+2\Big(2+n_{2}\Big)d_{G_1}(u)\Big]+\sum_{uv\in E(G_2)}\Big[\Big(d_{G_2}(u)\\&+d_{G_2}(v)\Big)^2+4m_{1}^2+4m_{1}\Big(d_{G_2}(u)+d_{G_2}(v)\Big)\Big]+\sum_{\substack{uv\in E(G_1\dotplus G_2)\\ u\in V_{s}(G_1)\\ v\in V(G_2)}}\Big[d_{G_2}(v)^2+\\&\Big(m_{1}+n_{2}+2\Big)^2+2\Big(m_{1}+n_{2}+2\Big)d_{G_2}(v)\Big]\\
			&=\sum_{\substack{uv\in E(S(G_1))\\ u\in V(G_1)\\ v\in V_{s}(G_1)}}d_{G_1}(u)^2+\Big(2+n_{2}\Big)^2\sum_{uv\in E(S(G_1))}1+2\Big(2+n_{2}\Big)\sum_{\substack{uv\in E(S(G_1))\\ u\in V(G_1)\\ v\in V_{s}(G_1)}}d_{G_1}(u)+\\&\sum_{uv\in E(G_2)}\Big(d_{G_2}(u)+d_{G_2}(v)\Big)^2+4m_{1}^2\sum_{uv\in E(G_2)}1+4m_{1}\sum_{uv\in E(G_2)}\Big(d_{G_2}(u)+d_{G_2}(v)\Big)\\&+\sum_{u\in V_{s}(G_1)}1\sum_{v\in V(G_2)}d_{G_2}(v)^2+\Big(m_{1}+n_{2}+2\Big)^2\sum_{u\in V_{s}(G_1)}1\sum_{v\in V(G_2)}1+2\Big(m_{1}+n_{2}+2\Big)\\&\sum_{u\in V_{s}(G_1)}1\sum_{v\in V(G_2)}d_{G_2}(v)\\
			&=\sum_{\substack{uv\in E(S(G_1))\\ u\in V(G_1)\\ v\in V_{s}(G_1)}}\Big[d_{G_1}(u)^2+\Big(4+2n_{2}\Big)d_{G_1}(u)\Big]+m_{1}^{\prime}\Big(2+n_{2}\Big)^2+HM_{1}(G_2)+4m_{1}^2m_{2}\\&+4m_{1}M_{1}(G_2)+m_{1}M_{1}(G_2)+m_{1}n_{2}\Big(m_{1}+n_{2}+2\Big)^2+4m_{1}m_{2}\Big(m_{1}+n_{2}+2\Big)\\
			&=m_{1}\Big[5M_{1}(G_2)+4m_{1}m_{2}\Big]+HM_{1}(G_2)+m_{1}^{\prime}\Big(2+n_{2}\Big)^2+m_{1}\Big(m_{1}+n_{2}+2\Big)\\&\Big[n_{2}\Big(m_{1}+n_{2}+2\Big)+4m_{2}\Big]+\sum_{uv\in E(S(G_1))}\Big[d_{G_1}(u)^2+\Big(4+2n_{2}\Big)d_{G_1}(u)\Big]
		\end{align*}
		
	which is our required result.
	\end{proof}
\end{thm}

Here is the example. By using the statement of theorem \ref{H2}, we get

\begin{align*}
	HM_{1}(P_{l}\dotplus C_{m})=&(l-1)\Big[20m+4m(l-1)\Big]+16m+2(l-1)(2+m)^2\\&+(l-1)\Big[(l-1)+m+2\Big]\Big[m\Big((l-1)+m+2\Big)+4m\Big]\\&+\sum_{uv\in E(S(P_l))}\Big[d_{P_l}(u)^2+(4+2m)d_{P_l}(u)\Big]
\end{align*}

\begin{thm} \label{R1}
	Let we have two simple connected graphs $G_{1}=(n_{1},m_{1})$ and $G_{2}=(n_{2},m_{2})$, then the reduced second Zagreb index of Corona join product $G_{1}\oplus G_{2}$ is given as
	\begin{align*}
		&RM_2(G_{1}\oplus G_{2})=RM_{2}(G_1)+n_{1}RM_{2}(G_2)+n_{1}\Big[n_{2}M_{1}(G_1)+n_{1}M_{1}(G_2)\Big]+\\&m_{1}n_{1}n_{2}\Big(n_{1}n_{2}-2\Big)+m_{2}n_{1}^2\Big(n_{1}-2\Big)+2n_{1}\Big[2m_{1}m_{2}+m_{1}n_{2}\Big(n_{1}-1\Big)+\\&m_{2}n_{1}\Big(n_{1}n_{2}-1\Big)\Big]+n_{1}^2n_{2}\Big[n_{1}^2n_{2}-n_{1}n_{2}-n_{1}+1\Big]
	\end{align*}
	
	\begin{proof}
		From the definition of the reduced second Zagreb index we have
		
		\[RM_2(G_{1}\oplus G_{2})=\sum_{uv \in E(G_{1}\oplus G_{2})}\Big[d_{G_{1}\oplus G_{2}}(u)-1\Big]\Big[d_{G_{1}\oplus G_{2}}(v)-1\Big]\]
		
		Now we apply the Lemma \ref{cjp}
		\begin{align*}
			&=\sum_{uv\in E(G_{1})}\Big(d_{G_{1}}(u)+n_{1}n_{2}-1\Big)\Big(d_{G_{1}}(v)+n_{1}n_{2}-1\Big)+n_{1}\sum_{uv\in E(G_{2})}\\&\Big(d_{G_{2}}(u)+n_{1}-1\Big)\Big(d_{G_{2}}(v)+n_{1}-1\Big)+n_{1}\sum_{\substack{uv\in E(G_{1}\oplus G_{2})\\ u\in V(G_1)\\ v\in V(G_2)}}\Big(d_{G_{1}}(u)+n_{1}n_{2}-1\Big)\\&\Big(d_{G_{2}}(v)+n_{1}-1\Big)\\
			&=\sum_{uv\in E(G_{1})}\Big[\Big[d_{G_{1}}(u)d_{G_{1}}(v)-d_{G_{1}}(u)-d_{G_{1}}(v)+1\Big]+n_{1}n_{2}\Big(d_{G_{1}}(u)+d_{G_{1}}(v)\Big)\\&+(n_{1}n_{2})^2-2n_{1}n_{2}\Big]+n_{1}\sum_{uv\in E(G_{2})}\Big[\Big[d_{G_{2}}(u)d_{G_{2}}(v)-d_{G_{2}}(u)-d_{G_{2}}(v)+1\Big]
		\end{align*}
		\begin{align*}
			\\&+n_{1}\Big(d_{G_{2}}(u)+d_{G_{2}}(v)\Big)+n_{1}^2-2n_{1}\Big]+n_{1}\sum_{\substack{uv\in E(G_{1}\oplus G_{2})\\ u\in V(G_1)\\ v\in V(G_2)}}\Big[\Big[d_{G_{1}}(u)d_{G_{2}}(v)+\\&\Big(n_{1}-1\Big)d_{G_{1}}(u)+\Big(n_{1}n_{2}-1\Big)d_{G_{2}}(v)+n_{1}^2n_{2}-n_{1}n_{2}-n_{1}+1\Big]\Big]\\
			&=\sum_{uv\in E(G_{1})}\Big[\Big[\Big(d_{G_{1}}(u)-1\Big)\Big(d_{G_{1}}(v)-1\Big)\Big]+n_{1}n_{2}\Big(d_{G_{1}}(u)+d_{G_{1}}(v)\Big)+(n_{1}n_{2})^2\\&-2n_{1}n_{2}\Big]+n_{1}\sum_{uv\in E(G_{2})}\Big[\Big[\Big(d_{G_{2}}(u)-1\Big)\Big(d_{G_{2}}(v)-1\Big)\Big]+n_{1}\Big(d_{G_{2}}(u)+d_{G_{2}}(v)\Big)\\&+n_{1}^2-2n_{1}\Big]+n_{1}\sum_{\substack{uv\in E(G_{1}\oplus G_{2})\\ u\in V(G_1)\\ v\in V(G_2)}}\Big[\Big[d_{G_{1}}(u)d_{G_{2}}(v)+\Big(n_{1}-1\Big)d_{G_{1}}(u)+\Big(n_{1}n_{2}-1\Big)d_{G_{2}}(v)\\&+n_{1}^2n_{2}-n_{1}n_{2}-n_{1}+1\Big]\Big]\\
			&=\sum_{uv\in E(G_{1})}\Big(d_{G_{1}}(u)-1\Big)\Big(d_{G_{1}}(v)-1\Big)+n_{1}n_{2}\sum_{uv\in E(G_{1})}\Big(d_{G_{1}}(u)+d_{G_{1}}(v)\Big)+\\&(n_{1}n_{2})^2\sum_{uv\in E(G_{1})}1-2n_{1}n_{2}\sum_{uv\in E(G_{1})}1+n_{1}\sum_{uv\in E(G_{2})}\Big(d_{G_{2}}(u)-1\Big)\Big(d_{G_{2}}(v)-1\Big)\\&+n_{1}^2\sum_{uv\in E(G_{2})}\Big(d_{G_{2}}(u)+d_{G_{2}}(v)\Big)+n_{1}^3\sum_{uv\in E(G_{2})}1-2n_{1}^2\sum_{uv\in E(G_{2})}1+n_{1}\sum_{u\in V(G_1)}d_{G_1}(u)\\&\sum_{v\in V(G_2)}d_{G_2}(v)+n_{1}\Big(n_{1}-1\Big)\sum_{u\in V(G_1)}d_{G_1}(u)\sum_{v\in V(G_2)}1+n_{1}\Big(n_{1}n_{2}-1\Big)\sum_{v\in V(G_2)}d_{G_2}(v)\\&\sum_{u\in V(G_1)}1+n_{1}^3n_{2}\sum_{u\in V(G_1)}1\sum_{v\in V(G_2)}1-n_{1}^2n_{2}\sum_{u\in V(G_1)}1\sum_{v\in V(G_2)}1-n_{1}^2\sum_{u\in V(G_1)}1\sum_{v\in V(G_2)}1\\&+n_{1}\sum_{u\in V(G_1)}1\sum_{v\in V(G_2)}1\\
			&=RM_{2}(G_1)+n_{1}n_{2}M_{1}(G_1)+m_{1}\Big(n_{1}n_{2}\Big)^2-2m_{1}n_{1}n_{2}+n_{1}RM_{2}(G_2)+n_{1}^2M_{1}(G_2)\\&+m_{2}n_{1}^3-2m_{2}n_{1}^2+4m_{1}m_{2}n_{1}+2m_{1}n_{1}n_{2}\Big(n_{1}-1\Big)+2m_{2}n_{1}^2\Big(n_{1}n_{2}-1\Big)+n_{1}^4n_{2}^2-n_{1}^3n_{2}^2\\&-n_{1}^3n_{2}+n_{1}^2n_{2}\\
		\end{align*}
		\begin{align*}
			&=RM_{2}(G_1)+n_{1}RM_{2}(G_2)+n_{1}\Big[n_{2}M_{1}(G_1)+n_{1}M_{1}(G_2)\Big]+m_{1}n_{1}n_{2}\Big(n_{1}n_{2}-2\Big)\\&+m_{2}n_{1}^2\Big(n_{1}-2\Big)+2n_{1}\Big[2m_{1}m_{2}+m_{1}n_{2}\Big(n_{1}-1\Big)+m_{2}n_{1}\Big(n_{1}n_{2}-1\Big)\Big]\\&+n_{1}^2n_{2}\Big[n_{1}^2n_{2}-n_{1}n_{2}-n_{1}+1\Big]
		\end{align*}
		which is our required result.
	\end{proof}	
\end{thm}

Here is the example. By using the statement of theorem \ref{R1}, we get
\begin{align*}
	RM_{2}(P_{l}\oplus C_{m})=&(l-3)+lm+l\Big[m(4l-6)+4lm\Big]+lm(l-1)(lm-2)\\&+ml^2(l-2)+2l\Big[2m(l-1)+m(l-1)^2+ml(lm-1)\Big]\\&+l^2m\Big[l^2m-lm-l+1\Big]
\end{align*}

\begin{thm} \label{R2}
Let we have three simple connected graphs $G_{1}=(n_{1},m_{1})$, \\$G_{2}=(n_{2},m_{2})$ and $S(G_{1})=(n_{1}^{\prime},m_{1}^{\prime})$, then the reduced second Zagreb index of sub-division vertex join $G_{1}\dotplus G_{2}$ is given as
	\begin{align*}
	&RM_{2}(G_{1}\dotplus G_{2})=RM_{2}(G_2)+m_{1}\Big[M_{1}(G_2)+m_{1}m_{2}-2m_{2}\Big]-m_{1}^{\prime}\Big(n_{2}+1\Big)\\&+m_{1}\Big(n_{2}+1\Big)\Big[2m_{2}+m_{1}n_{2}-n_{2}\Big]+\Big(n_{2}+1\Big)\sum_{uv\in E(S(G_1))}d_{G_1}(u)
	\end{align*}

\begin{proof}
		From the definition of the reduced second Zagreb index we have
	
	\[RM_2(G_{1}\dotplus G_{2})=\sum_{uv \in E(G_{1}\dotplus G_{2})}\Big[d_{G_{1}\dotplus G_{2}}(u)-1\Big]\Big[d_{G_{1}\dotplus G_{2}}(v)-1\Big]\]
	Now we apply the Lemma \ref{svjg}
	\begin{align*}
		&=\sum_{\substack{uv\in E(S(G_1))\\ u\in V(G_1)\\ v\in V_{s}(G_1)}}\Big[\Big(d_{G_1}(u)-1\Big)\Big(2+n_{2}\Big)\Big]+\sum_{uv\in E(G_2)}\Big[\Big(d_{G_2}(u)+m_{1}-1\Big)\\&\Big(d_{G_2}(v)+m_{1}-1\Big)\Big]+\sum_{\substack{uv\in E(G_1\dotplus G_2)\\ u\in V_{s}(G_1)\\ v\in V(G_2)}}\Big[\Big(2+n_{2}-1\Big)\Big(d_{G_2}(v)+m_{1}-1\Big)\Big]\\
	\end{align*}
	\begin{align*}
	&=\Big(n_{2}+1\Big)\sum_{\substack{uv\in E(S(G_1))\\ u\in V(G_1)\\ v\in V_{s}(G_1)}}\Big(d_{G_1}(u)-1\Big)+\sum_{uv\in E(G_2)}\Big[\Big[d_{G_2}(u)d_{G_2}(v)-d_{G_2}(u)\\&-d_{G_2}(v)+1\Big]+m_{1}\Big(d_{G_2}(u)+d_{G_2}(v)\Big)+m_{1}^2-2m_{1}\Big]+\sum_{\substack{uv\in E(G_1\dotplus G_2)\\ u\in V_{s}(G_1)\\ v\in V(G_2)}}\Big[\Big(n_{2}+1\Big)d_{G_2}(v)\\&+m_{1}\Big(n_{2}+1\Big)-\Big(n_{2}+1\Big)\Big]\\
	&=\Big(n_{2}+1\Big)\sum_{\substack{uv\in E(S(G_1))\\ u\in V(G_1)\\ v\in V_{s}(G_1)}}d_{G_1}(u)-\Big(n_{2}+1\Big)\sum_{\substack{uv\in E(S(G_1))}}1+\sum_{uv\in E(G_2)}\Big(d_{G_2}(u)-1\Big)\\&\Big(d_{G_2}(v)-1\Big)+m_{1}\sum_{uv\in E(G_2)}\Big(d_{G_2}(u)+d_{G_2}(v)\Big)+m_{1}^2\sum_{uv\in E(G_2)}1-2m_{1}\sum_{uv\in E(G_2)}1+\\&\Big(n_{2}+1\Big)\sum_{v\in V(G_2)}d_{G_2}(v)\sum_{u\in V_{s}(G_1)}1+m_{1}\Big(n_{2}+1\Big)\sum_{v\in V(G_2)}1\sum_{u\in V_{s}(G_1)}1-\Big(n_{2}+1\Big)\\&\sum_{v\in V(G_2)}1\sum_{u\in V_{s}(G_1)}1\\
	&=\Big(n_{2}+1\Big)\sum_{\substack{uv\in E(S(G_1))\\ u\in V(G_1)\\ v\in V_{s}(G_1)}}d_{G_1}(u)-m_{1}^{\prime}\Big(n_{2}+1\Big)+RM_{2}(G_2)+m_{1}M_{1}(G_2)+m_{1}^2m_{2}\\&-2m_{1}m_{2}+2m_{1}m_{2}\Big(n_{2}+1\Big)+m_{1}^2n_{2}\Big(n_{2}+1\Big)-m_{1}n_{2}\Big(n_{2}+1\Big)\\
	&=RM_{2}(G_2)+m_{1}\Big[M_{1}(G_2)+m_{1}m_{2}-2m_{2}\Big]-m_{1}^{\prime}\Big(n_{2}+1\Big)+m_{1}\Big(n_{2}+1\Big)\\&\Big[2m_{2}+m_{1}n_{2}-n_{2}\Big]+\Big(n_{2}+1\Big)\sum_{uv\in E(S(G_1))}d_{G_1}(u)
	\end{align*}
    
which is our required result.
 \end{proof}
\end{thm}

Here is the example. By using the statement of theorem \ref{R2}, we get
\begin{align*}
	&RM_{2}(P_{l}\dotplus C_{m})=m+(l-1)\Big[4m+m(l-1)-2m\Big]-2(l-1)(m+1)\\&+(l-1)(m+1)\Big[2m+m(l-1)-m\Big]+(m+1)\sum_{uv\in E(S(P_l))}d_{P_l}(u)
\end{align*}

\section{Conclusion}

In this paper we proposed the two new variants of special graph products and then find their exact expressions for the forgotten topological index, first hyper Zagreb index and reduced second Zagreb index. In future other degree and distance based topological indices of these graph operations can be found. The obtained results may help to construct and investigate the topological indices of complex networks structures.

\end{document}